\documentclass[runningheads]{llncs}
\usepackage{graphicx}
\usepackage{amssymb,amsmath}
\usepackage{algorithm}
\usepackage{algpseudocode}
\usepackage[noshellescape]{abc}
\begin{document}
\title{Melody and Variation Generation Through KAM Theory}
%
%
\author{Octavio A. Agustín-Aquino\inst{1}\orcidID{0000-0002-0556-6236} \and Alicia Santiago-Santos\inst{1}\orcidID{0000-0002-1888-7135}}
%
\authorrunning{O. A. Agustín-Aquino and A. Santiago-Santos}
%
\institute{Instituto de Física y Matemáticas, Universidad Tecnológica de la Mixteca, \url{http://www.utm.mx}}
%
\maketitle              
\begin{abstract}
We use the Newton iteration of KAM theory for diffeomorphisms of the circle as a source for melody generation and its variations.

\keywords{Music \and melody generation \and discrete dynamical systems \and KAM theory.}
\end{abstract}
\section{Introduction}

Dynamical systems can be used to generate melodies or its variations, as was shown by the pioneering work of Pressing \cite{jP88} for certain discrete and discretized systems and Dabby \cite{dD96} for the case of altering the tonal characteristics of a piece by chaotic trajectories; see \cite{kaliakatsos2013chaos} for a survey on this and related topics.

In a nutshell, the idea is to associate a sequence of pitches $p_{k}$ to a sequence of points $x_{k}$ in a space $X$ such that
\[
 x_{k+1} = f(x_{k})
\]
for a given $f:X\to X$ defined by a dynamical system. It is undesirable for the behaviour of the sequence $\{x_{k}\}_{k\in\mathbb{N}\cup\{0\}}$ to be too ``predictable'', since it would be uninteresting. On the contrary, if it is too ``chaotic'' then it would yield melodies indistinguishable from ``gibberish'', so we would like to have some control between these two extremes.

Here we propose to use KAM theory\footnote{KAM stands for Kolmogorov, Arnold and Moser, who developed it.} for the diffeomorphisms of the circle \cite{Arnold2009,ref_proc1} as a reasonable alternative not attested in existing literature to accomplish this control, by perturbing an original melody constructed with a simple dynamical system and gradually returning to it. Since the behaviour of the system and its perturbation depends on two easily tunable parameters and the choice of one function, we deem it as a good candidate for melody and variation generation.

In Section 2 we establish some notation and terminology that we will use throughout this paper. In Section 3 Arnold's theorem is stated after certain preliminaries. This allows us to describe the variation algorithm in Section 4 and provide some notes regarding its implementation in Octave. A few examples are given in Section 5.

Finally, in Section 6, some conclusions and a glimpse of future steps of this work are reached.

\section{Some basic notions}

In this section we introduce some general notation and the concept of dynamical system, as well as certain terminology we will use frequently. For further details concerning circle mappings or dynamical systems, see \cite{FG22} and \cite{ref_Holmgren}, respectively.

\begin{definition}
Given a topological space $X$ and a semigroup $(G,*)$, a \emph{dynamical system} is a continuous function $\psi:G\times X\to X$ such that
\begin{enumerate}
\item for every $x\in X$ we have $\psi(0,x)=x$, where $0$ is the neutral element of $G$ and
\item for every $t,s\in G$ and $x\in X$ we have $\psi(t,\psi(s,x))=\psi(t*x,x)$.
\end{enumerate}
\end{definition}

In particular, we will restrict our attention to \emph{discrete} dynamical systems with $G=(\mathbb{N}\cup \{0\},+)$ and $\psi$ defined via a continuous function $f:X\to X$ by
\[
 \psi(k,x) = f^{k}(x),
\]
where $f^{k}$ is the $k$-fold composition of $f$ with itself \cite[Definition 4.5, p. 33]{ref_Holmgren}. Given $x\in X$, we define the its \emph{orbit} under $f$ as the set
\[
\mathcal{O}(x,f) = \{f^{k}(x):k\in \mathbb{N}\cup\{0\}\}.
\]

Now we specify $f$ further, taking it as a diffeomorphism of the circle $f:S^{1}\to S^{1}$. There exists a homeomorphism $\phi:\mathbb{R}\to\mathbb{R}$ such that
\[
f(\exp(2\pi i t))=\exp(2\pi i \phi(t))
\]
for every $t\in\mathbb{R}$, which is called a \emph{lift} of $f$. 

A fundamental example of a diffeomorphism of the circle is a \emph{rotation}. The lift of a circle rotation $r_{\alpha}$ by an angle $\alpha$ will be denoted by $ R_{\alpha}$. Thus,
\[
 R_{\alpha}(x) = x + \alpha
\]
since, as a matter of fact, for any diffeomorphism of the circle $f$ with lift $\phi$ the limit
\[
\rho(f) = \lim_{n\to\infty}\frac{\phi ^{n}(x)-x}{n}
\]
exists and it is called the \emph{rotation number} of $f$. As we will see soon, any diffeomorphism of the circlr $f$ is essentially equivalent to $r_{\rho(f)}$.

Moreover, the behaviour of the orbits under a rotation $r_{\alpha}$ is completely understood: if $\alpha\in \mathbb{Q}$ then $\mathcal{O}(x,r_{\alpha})$ is a finite set; otherwise, it is dense in the circle. It is a good moment to mention that an irrational number can be approximated arbitrarily well by some rational number $\frac{m}{n}$, but in general such an approximation cannot get closer than a distance proportional to $\frac{1}{n^{2}}$. For our purposes we need some quantification on this fact.

\begin{definition}\label{Def:AN} A number $\rho\in\mathbb{R}$ is of type $(K,\nu)$ if there exist positive numbers $K$ and $\nu$ such that
\[
 \left|\rho-\frac{m}{n}\right|>\frac{K}{|n|^{\nu}}
\]
for al pairs of integers $(m,n)$.
\end{definition}

The number $\nu$ from Definition \ref{Def:AN} is called the \emph{irrationality measure} of $\rho$ \cite[p. 28]{yB04}. Roth proved that for any irrational algebraic number its irrationality measure is $2$ \cite[p. 28, Theorem 2.1]{yB04}; thus, for example, $\varphi = \frac{\sqrt{5}+1}{2}$ is of type $(K,2)$. In \cite{jS06} there is a nice proof that $e$ has irrationality measure $2$ as well. The irrationality measure of $\pi$ is not known; as of this writing the best current upper bound on it is approximately $7.103205334$ \cite{ZZ20}.

\section{Arnold's theorem}

In this section we deal first with the preliminaries to understand Arnold's theorem statement and explain what the \emph{Newton iteration} is in this context.

Let us begin considering diffeomorphisms of the circle with lifts of the form
\[
\phi(x) = x + \alpha + \eta(x),
\]
where $\alpha$ is an irrational number, $\eta:\mathbb{R}\to\mathbb{R}$ is periodic with period $1$ and $\frac{d}{dx}\eta(x)>-1$ (so orientation is preserved). Furthermore, we will suppose
\[
\sup_{|\mathrm{Im}(z)|<\sigma}|\eta(z)| \leq \infty.
\]

We are looking for a change of variables $x=H(\xi)$ such that
\[
 H^{-1}\circ \phi \circ H = R_{\alpha},
\]
or, equivalently,
\begin{equation}\label{Eq:Cambio}
 \phi \circ H = H\circ R_{\alpha}.
\end{equation}

If we demand from $H$ just to be an homeomorphism, then Denjoy's theorem \cite[Theorem 3.4, p. 48]{FG22} suffices to guarantee its existence. But, if we also require $H$ to be analytical, then Arnold's theorem not only tells us that it exists (under certain relatively mild conditions) but we can also profit from its proof to successively build approximations to $H$.

Let us elaborate: Suppose $H$ is of the form $H(x) = x+h(x)$. Substituting this in \eqref{Eq:Cambio} leads to
\begin{align*}
 \phi(H(x))=\phi(x+h(x)) &= x+h(x)+\alpha+\eta(x) \\
 &=H(R_{\alpha}(x)) = x+\alpha + h(x+\alpha),
\end{align*}
thus
\[
h(x+\alpha)-h(x) = \eta(x).
\]

This implies that the Fourier coefficients of $\hat{\eta}(n)$ satisfy
\[
 \hat{h}(n) = \frac{\hat{\eta}(n)}{e^{2\pi i n \rho}-1}
\]
for $n\in \mathbb{Z}\setminus\{0\}$. We can now write
\[
h(x) = \sum_{n\in\mathbb{Z},n\neq 0}\frac{\hat{\eta}(n)}{e^{2\pi i n \rho}-1}e^{-2\pi i n x}.
\]

Unless $\hat{\eta}(0) = 0$ this does not completely solve the problem of finding $H$, but the approximate $h$ obtained this way suffices, for we can now define 
\[
\phi_{1}(x) = H^{-1}\circ \phi\circ H
\]
so we can calculate an $h_{1}(x)$ through the Fourier coefficients of $\eta_{1}(x) = \phi_{1}(x)-\alpha-x$, and continue inductively with $H_{1}(x)=x+h_{1}(x)$ and
\[
\phi_{n}(x) = H_{n-1}^{-1}\circ \phi_{n-1}\circ H_{n-1}
\]
for $n\geq 2$. This is the so-called \emph{Newton iteration} for the calculation of the change of coordinates. Now we can state Arnold's theorem.

\begin{theorem}[Arnold, 1965, {\cite[Theorem 2]{Arnold2009}}] If
\begin{enumerate}
\item the number $\rho$ is of type $(K,\nu)$,
\item the diffeomorphism
\[
 \phi(x) = x+\rho+\eta(x)
\]
has rotation number $\rho$, and
\item for some $\epsilon(K,\nu,\sigma)$ we have $\sup_{|\mathrm{Im}(z)|<\sigma}|\eta(z)| \leq \epsilon(K,\nu,\sigma)$,
\end{enumerate}
then there is an $N$ such that the change of variables $H=H_{1}\circ\cdots\circ H_{N}$ is analytical and conjugates $\phi$ to $R_{\rho}$.
\end{theorem}

\section{Algorithm and implementation details}

Our approach to generate a sequence of pitches is to take $x_{0}=p_{0}=0$ and then, given a circle diffeomorphism with lift $\phi:\mathbb{R}\to\mathbb{R}$, to calculate
\[
p_{k} = \lfloor 12\bmod(\phi^{k}(x_{0}),1)\rfloor,
\]
thus we obtain a melody where each pitch is codified as semitones from the tonic. It is to be noted that, since in the long run $f$ behaves like a rotation when it comes to the orbit, this tends to stabilize in a pattern that shifts slowly, because of the coarse grain of the $12$-tone scale. Of course it is possible to consider other scales that subdivide the octave with a different or finer grain, but we will not pursue such a direction in this paper.

Our diffeomorphism of departure is one with a lift among Arnold's family \cite[p. 248]{AP90}
\begin{align*}
\phi :\mathbb{R}&\to\mathbb{R},\\
x&\mapsto x+\alpha+\frac{\epsilon}{2\pi}\sin(2\pi x),
\end{align*}
where $\epsilon$ is chosen freely (but small, say, $\epsilon<1$). The detailed procedure is contained in Algorithm \ref{alg:KAM}. It was implemented\footnote{See the authors' GitHub directory.} in Octave 5.2 as a function of $\alpha$, $\epsilon$, $n$ and $M$, and of course it is possible to modify the number of samples of step 8. The estimation of Fourier coefficients of step 9 is done with a right-point Riemann sum and the solution step 14 is done with the native Octave function \texttt{fsolve}, with the initial estimation chosen as $x_{k-1}+\alpha$. The output is the pure numerical result, before the codification as semitones from the tonic is done.

\begin{algorithm}
\caption{Melody variation through KAM theory}\label{alg:KAM}
\begin{algorithmic}[1]
\Require $\alpha,\epsilon \in \mathbb{R}_{+}$ with $0<\epsilon< 1$, $n,M\in\mathbb{N}$, where $n$ is the length of the melody and $M$ is the number of variations.
\Ensure $M$ sets $\{p_{k}\}$ of cardinality $n$ that correspond to pitches of $M$ variations of a melody.
\State $\phi \gets x+\alpha+\frac{\epsilon}{2\pi}\sin(2\pi x)$.
\State $x_{0},p_{0}\gets 0$
\For{$1\leq k \leq n-1$}
 \State $x_{k} \gets \phi(x_{k-1})$.
 \State $p_{k} \gets \lfloor 12\bmod(x_{k},1)\rfloor$.
\EndFor
\For{$1\leq s \leq M$}
 \State Sample the function $\eta(x)=\phi(x)-\alpha-x$ in $[0,1]$.
 \State Use the sample to estimate the Fourier coefficients $\hat{\eta}(j)$ for $1\leq j\leq N$.
 \State $h\gets \sum_{-N\leq j\leq N, j\neq 0} \frac{\hat{\eta}(j)}{e^{2\pi i j \alpha}-1}e^{-2\pi i j x}$.
 \State $H(x)\gets x+h(x)$.
 \For{$1\leq k \leq n-1$}
  \State $u\gets \phi\circ H (x_{k-1})$.
  \State Solve $H(x_{k}) = u$ for $x_{k}$.
  \State $p_{k} \gets \lfloor 12\bmod(x_{k},1)\rfloor$.
 \EndFor
 \State \Return $\{p_{k}\}$.
 \State $\phi\gets H^{-1}\circ \phi \circ H$.
\EndFor
\end{algorithmic}
\end{algorithm}

\section{Some examples}

The code to obtain the results of the following three examples can be found in the authors' GitHub directory.

\begin{example}
If we choose $\alpha = e$ and $\epsilon = 0.5$ in Arnold's diffeomorphism, the first ten pitches generated through $\phi_{1}^{k}$, $1\leq k \leq 9$, are
\[
 0, 8, 4, 1, 11, 7, 3, 0, 9, 5
\]
which correspond to the melody seen in Fig. \ref{fig1} within the octave.

\begin{figure}
\includegraphics[width=\textwidth]{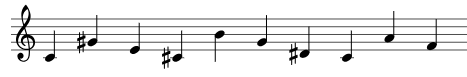}
\caption{First melody.} \label{fig1}
\end{figure}

We can calculate exactly the Fourier coefficients of $\eta_{1}$, and thus $H_{1}$, which is
\[
 H_{1}(x) = x - \frac{\epsilon}{2\pi}\frac{\cos(2\pi x-e\pi)}{2\sin(e\pi)}.
\]

By calculating with $\phi_{2}=H_{1}^{-1}\circ \phi_{1} \circ H_{1}$, we obtain the pitches
\[
 0, 8, 5, 2, 10, 7, 4, 1, 9, 6,
\]
namely the melody seen in Fig. \ref{fig2}.


\begin{figure}
\includegraphics[width=\textwidth]{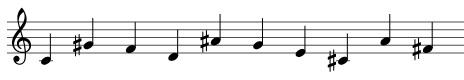}
\caption{Second melody, obtained as a variation of the first.} \label{fig2}
\end{figure}

It is also possible to calculate $H_{2}$ exactly \cite[p. 214]{Arnold2009}, which is
\[
H_{2}(x) = x-\frac{\epsilon}{2\pi}\frac{\cos(2\pi x-e\pi)}{2\sin(e\pi)}
+\frac{\epsilon^{2}}{4\pi^{2}}\frac{\sin(4\pi x-e\pi)}{4\sin(e\pi)\sin(2e\pi)},
\]
but the resulting change with respect to the second variation is nil.
\end{example}

\begin{example}
If we use the same diffeomorphism of the previous example but we increase the perturbation to $\epsilon = 0{.}8$, then we obtain the melodies seen in Figure \ref{fig3} generated by the exact change of coordinates of two iterations, in this case separated by barlines.

\begin{figure}
\includegraphics[width=\textwidth]{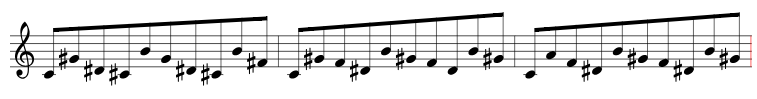}
\caption{Three melodies (separated by barlines) obtained taking $\epsilon = 0{.}8$, $\alpha = e$, $n=10$ and $M=3$ in Algorithm \ref{alg:KAM}.} \label{fig3}
\end{figure}
\end{example}

\begin{example}
If we run the implementation of Algorithm \ref{alg:KAM} with $\alpha = \pi$, $\epsilon = 0{.}8$, $n=10$ and $M=7$ with $64$ equally spaced point samples and keeping $10$ Fourier coefficients, then we obtain the following output after the semitone codification, displayed as a matrix,
\[
P = \begin{pmatrix}
    0   & 1  &  4  &  7 & 8 & 8 & 8 & 8 & 9 & 9\\
  0  &  0 &   1  &  2 & 4 & 5 & 6 & 6 & 6 & 7\\
 0   & 2  &  3  &  5 & 7 & 6 & 8 & 8 & 5 & 6\\
  0   & 1  &  3  &  5 & 6 & 8 & 6 & 8 & 9 &  11\\
  0   & 1  &  3  &  5 & 6 & 8 &  10 & 0 & 2 & 3\\
  0   & 1  &  3  &  5 & 6 & 8 &  10 &  11 & 1 & 3\\
  0   & 1  &  3  &  5 & 6 & 8 &  10 &  11 & 1 & 3\\
\end{pmatrix}.
\]

In this case the convergence towards the rotation is patent. 
\end{example}

\section{Conclusions and future work}

As we can readily see the variation obtained by applying a change of coordinates to a circle diffeomorphism in the Arnold family is noticeable yet subtle, similar in spirit to Glass' ``repetitive structures'' \cite[Chapter ``Paris'']{pG15}. The numerical implementation introduces additional relatively small ``instabilities'', but this could be seen not as a flaw but as a feature!

With the restriction to this type of dynamical system we can only control one aspect of melodic variation (pitch, in this case), but KAM theory also works for integrable Hamiltonian systems, where in principle we can also obtain control over note length or (musical) dynamics, for example. The research in this direction is ongoing.

%
%
%
\bibliographystyle{splncs04}
\bibliography{mybibliography}

\end{document}